\documentclass[11pt]{article}
\usepackage{amsmath,amssymb,mathrsfs}
\date{}

\def\tr{\mathop{\rm tr}}

\renewcommand{\mathcal}{\mathscr}
\def\cl#1{{\mathscr #1}}
\newcommand{\dlines}{\displaylines}
\def\ep{\varepsilon}

\def\paref#1{(\ref{#1})}

\newcommand{\field}[1]{\mathbb{#1}}
\newcommand{\R}{\field{R}}
\newcommand{\qed}{\par\hfill$\blacksquare$\par\noindent\ignorespaces}

\newcommand{\RR}{\field{R}}

\let\pa=\partial

\newcommand{\e}{{\rm e}}

\def\tin#1{\par\noindent\hskip3em\llap{#1\enspace}\ignorespaces}
\def\cl#1{{\mathcal #1}}

\def\P{{\rm P}}

\def\tfrac#1#2{{\textstyle\frac {#1}{#2}}}
\newtheorem{theorem}{Theorem}

\numberwithin{equation}{section}

\newtheorem{assum}{Assumption}

\newtheorem{ex}[theorem]{Example}

\newtheorem{prop}[theorem]{Proposition}
\newtheorem{remark}[theorem]{Remark}

\begin{document}
\title{\Huge\sc On Sharp Large Deviations for the bridge of a general Diffusion}
\author{Paolo Baldi, Lucia Caramellino, Maurizia Rossi
\\
{\it Dipartimento di Matematica, Universit\`a
di Roma Tor Vergata, Italy}\\
} \date{}\maketitle

\abstract{
We provide sharp Large Deviation estimates for the probability of exit from a domain for the  bridge of a
$d$-dimensional general diffusion process $X$, as the conditioning time tends to $0$.
This kind of results is motivated by applications to numerical simulation.
In particular we investigate the influence of the drift $b$ of $X$. It turns out that
the sharp asymptotics for the exit time probability are independent of the drift,
provided $b$ enjoyes a simple condition that is
always satisfied in dimension $1$. On the other hand, we show that
the drift can be influential
if this assumption is not satisfied.
}
\smallskip

\smallskip


\section{Introduction}\label{intro}
Even if the subject of Large Deviations was not one of the most visited among the many objects of investigation in the large scientific production of Marc Yor, he was able to provide three original contributions in this field (\cite{MR2063378}, \cite{yor-arcsine}, \cite{yor-clock}). On the other hand bridges and conditioned processes have been at the heart of many of his most important contributions.
In this short note we investigate some points concerning the asymptotics of conditioned processes when the conditioning time goes to $0$.

The investigation of Large Deviation and sharp Large Deviation estimates in this context goes back to \cite{Bal:95}, where it was motivated by applications to simulation. This line of research was continued in the subsequent years
(\cite{MR1849251},\cite{gobet},\cite{BC2002-MR1925452}) but when needed by the numerical application usually the conditioned Diffusion
was approximated by the bridge of the Diffusion with its coefficients frozen, thus taking advantage of the well
known asymptotics of the Brownian bridge. The need of more accurate estimates is now prompted by applications to finance, mainly in connection with
stochastic volatility models.

Actually when simulating the path of a stochastic process (with the Euler scheme e.g.) which is to be killed at the exit from some domain $D$, it is important to be able to compute the probability for the conditioned Diffusion with $X_{t_n}=x$, $X_{t_{n+1}}=y$ to exit from $D$ in the time interval $[t_n,t_{n+1}]$, where $t_n,t_{n+1}$ are consecutive times in the time grid and $X_{t_n}, X_{t_{n+1}}$ denote the corresponding simulated positions (see \cite{Bal:95}, p.1645-1646 for a more complete explanation).

In the remainder of this paper most of the time, when speaking of {\it asymptotics} we shall mean in the sense of the
conditioning time to go to $0$ and the terms {\it bridge} or {\it conditioned Diffusion} shall mean the same thing.

To find a sharp estimate for the exit from a given domain of the bridge of a general Diffusion process is a goal that requires
some work. In this note we wish to investigate a minor point in this direction.
It has been proved (\cite{bailleul}, \cite{BCR2}) that the (non sharp) Large Deviation asymptotics
for conditioned Diffusions do not depend on the drift  $b$ of the non conditioned process  $X$ as
in \paref{SDE1}. It has been a general belief that this remains true also for the sharp asymptotics of
the bridge of a Diffusion. This was actually proved for a large family of
one dimensional Diffusion processes in \cite{BC2002-MR1925452}. Our object is to prove that this
property actually holds in general
in the multidimensional situation, provided the drift satisfies a simple condition, that is
of course always satisfied in dimension $1$.

Our results stem from two main tools: the asymptotics of W.H. Fleming and M.R. James
\cite{FleJ:92} and the asymptotics for small time of the transition density of a diffusion process of S. A.
Mol{\v{c}}anov \cite{molchanov-MR0413289} together with those of G. Ben Arous \cite{benarous}.

Our goal here is mainly to put forward some ideas and techniques,
without trying to look for minimal regularity assumptions.

\section{Conditioned Diffusions}\label{CD}
Let $X$ be a $d$-dimensional (possibly inhomogeneous) diffusion process with transition density
$p$. The conditioned Diffusion
given $X_t=y$, $t>0$, is the one associated to the transition density
$$
\widehat
p(u,v,x,z)=\frac{p(u,v,x,z)p(v,t,z,y)}{p(u,t,x,y)}\ ,\quad 0\le u\le v<t,\  x,z\in\R^d\ .
$$
Remark that this is a time inhomogeneous transition density, even if $X$ was time homogeneous.
Let us assume moreover that $X$ is the solution of the Stochastic Differential Equation (SDE)
\begin{equation}\label{SDE1}
dX_t=b(X_t)\, dt+\sigma(X_t)\, dB_t
\end{equation}
(we consider therefore a process $X$ that is time homogeneous) and let us denote by $L$ its generator
$$
L=\frac 12\sum_{i,j=1}^d a_{ij}(z)\frac {\pa^2}{\pa z_i\pa z_j}+\sum_{i=1}^d
b_i(z)\frac {\pa}{\pa z_i}\ ,
$$
where, as usual, $a=\sigma\sigma^*$. By a straightforward computation, the generator $\widehat L_v$ of the conditioned Diffusion is
$$
\widehat L_v=L+\sum_{i=1}^d\widehat b_i^{y,t}(z,v)\frac{\pa\hfil}{\pa z_i}\ ,\quad 0\le v<t\ ,
$$
where
$$
\widehat b_i^{y,t}(z,v)=\frac1{p(t-v,z,y)}
\sum_{i,j=1}^d a_{ij}(z)\frac{\pa\hfil}{\pa z_j}p(t-v,z,y)=
\sum_{i,j=1}^d a_{ij}(z)\frac{\pa\hfil}{\pa z_j}\log p(t-v,z,y)\ ,
$$
i.e.
\begin{equation}\label{grad}
\widehat b^{y,t}(z,v)=a(z)\nabla_z p(t-v,z,y)\ .
\end{equation}
The conditioned Diffusion has therefore the same distribution as the solution of the SDE
$$
d\xi_v=\bigl(b(\xi_v)+\widehat b^{y,t}(\xi_v,v)\bigr)\, dv+\sigma(\xi_v)\,dB_v
$$
for $u<t$. Let
$\eta^t_v=\xi_{vt}$, the time changed conditioned Diffusion so that it is defined on $[0,1[$.
The new Diffusion $\eta^t$ is the solution of
\begin{equation}\label{time-changed}
d\eta^t_v=\bigl(b(\eta^t_v)+\widehat b^{y,t}(\eta^t_v,tv)\bigr)t\, dv+
\sqrt{t}\sigma(\eta^t_v)\,dB_v
\end{equation}
(with respect to a possibly different Brownian motion). We can therefore obtain estimates concerning the conditioned Diffusion
using the Ventsel-Freidlin Large Deviation estimates as soon as we are able to compute the limit
\begin{equation}\label{limite}
\widehat b^{y}(z,v):=\lim_{t\to 0}t\,\widehat b^{y,t}(z,tv)
\end{equation}
uniformly on compact sets (\cite{az-stflour,priouret,baldi-chaleyat}) and prove that the limit function $\widehat b^{y}$ is smooth enough.
Recall that
$$
\widehat b^{y,t}(z,tv)=\sum_{i,j=1}^d a_{ij}(z)\frac{\pa\hfil}{\pa z_j}\log p(t(1-v),z,y)=
a(z)\,\nabla_z \log p(t(1-v),z,y)\ .
$$
Let us denote $\widehat\P^{y,t}_{x,s}$ the law of $\eta^t$ with the starting condition $\eta^t_s=x$.
Let $D\subset\RR^d$ be an open set with a smooth boundary  and   $\tau=\tau_{x,s}$ the
exit time from $D$.  The Ventsel-Freidlin theory of  Large Deviations states that
\begin{equation}\label{ell}
\lim_{t\to0}t\log\widehat\P^{y,t}_{x,s}(\tau<1)=-\inf_{\gamma(s)=x,\tau(\gamma)<1}I_s(\gamma)=:-\ell_{x,s}\ ,
\end{equation}
 $I_s$ denoting the rate Ventsel-Freidlin function
\begin{equation}\label{I1}
I_s(\gamma)=\frac 12 \int_s^1
\langle{a(\gamma_v)^{-1}(\dot\gamma_v-\widehat b^{y}(\gamma_v,v))},{\dot\gamma_v-\widehat b^{y}(\gamma_v,
v)}\rangle \, dv\ ,
\end{equation}
where $a=\sigma\sigma^*$ as above.
\section{The sharp asymptotics}
We want to prove the stronger result
\begin{equation}\label{dev1}
q_t(x,s):=\widehat\P^{y,t}_{x,s}(\tau<1)\sim c_{x,s}\,\e^{-\ell_{x,s}/t}\ ,
\end{equation}
as $t\to 0$, for some constant $c_{x,s}> 0$.
We stress that it is relevant that $c_{x,s}$
is a constant independent of $t$ (see Remark \ref{relevant} for more comments).

We shall investigate the situation where the positions $x$ (the starting point of the process)
and $y$ (the conditioning position) are close to each other.
This is justified by the application mentioned in \S\ref{intro}:
$x$ and $y$ being consecutive positions in a simulation scheme it should be safe to
assume that they are  not far one from the other.

The computation of the asymptotics \paref{dev1} was performed in \cite{Bal:95} in the case where
$X$ is a multidimensional Brownian motion. The idea there was to take advantage of the results of W.A.~Fleming and M.~James
\cite{FleJ:92}.
Let us recall these estimates. Let $X^\ep$ be the solution of
\begin{equation}\label{be}
\begin{cases}
dX^\ep_t=\widetilde b_\ep(X^\ep_t,t)\,dt+\sqrt{\ep}\,\sigma(X^\ep_t)\,dB_t\ ,\quad t>s&\cr
X^\ep_s=x\in D\ .&\cr
\end{cases}\end{equation}
Let $T>0$ be fixed and let us assume that
$$
\lim_{\ep\to0}\widetilde b_\ep(x,s)=\widetilde b(x,s)
$$
uniformly for $(x,s)$ on the compact sets of $D\times [0,T]$.
Let us define the function $u:D\times [0,T[\to\RR$ by
\begin{equation}\label{eq-for-u}
u(x,s)=\inf_{\gamma(s)=x,\tau(\gamma)<T} \widetilde I_s(\gamma)
\end{equation}
where $\widetilde I_s$ is, similarly as in \paref{I1},
\begin{equation}\label{Itilde}
\widetilde I_s(\gamma)=\frac 12 \int_s^{T}
\langle{a(\gamma_t)^{-1}(\dot\gamma_t-\widetilde  b(\gamma_t,t))},{\dot\gamma_t-\widetilde  b(\gamma_t,
t)}\rangle \, dt\ .
\end{equation}
It can be shown that $u$ is the solution of the Hamilton-Jacobi problem
\begin{equation}\label{HJ1}
\begin{cases}
\displaystyle{\partial u\over\partial s}+\langle \widetilde b,\nabla u\rangle-\frac 12\langle a \nabla u,\nabla u\rangle=0& {\rm in\ }D\times ]0,T[\cr
u(x,s)=0&{\rm on\ }\partial D\times [0,T]\cr
u(x,s)\to +\infty&{\rm as\ }s\nearrow T,\  x\in D\cr
\end{cases}
\end{equation}
to be considered in the sense of viscosity solutions (\cite{FS})
and in the classical sense at each point of differentiability of $u$.

Now let $N\subset D\times [0,T'],\, T'<T$  and define
\begin{equation}\label{betaa}
\beta(x,s)=\widetilde b(x,s)-a(x)\nabla u(x,s),\qquad (x,s)\in \overline N\ .
\end{equation}
Let $\gamma_{x,s}$ be the solution of
\begin{equation}\label{e159}
\begin{cases}
\dot\gamma_{x,s}( t)=\beta(\gamma_{x,s}( t),t)&\cr
\gamma_{x,s}(s)=x\ &\cr
\end{cases}
\end{equation}
and set $t^*_{x,s}=\inf\{t>s, (\gamma_{x,s}( t),t)\not\in N\}$, moreover define
$$
\Gamma_1=\{(\gamma_{x,s}(t^*_{x,s}),t^*_{x,s}),(x,s)\in N\}\ .
$$
\begin{assum}\label{(A)} {\sl {a)} $N$ is an open set;
\tin {b)} $u\in\cl C^\infty(\overline N)$;
\tin{c)} $N$ is a  Region of Strong Regularity (RSR)
w.r.t. $\beta$, i.e.
$\Gamma_1$ is a $C^\infty$ manifold, relatively open in $\pa N$,
$(\gamma_{x,s}(t),t)_{t\in [s, t^*_{x,s}]}$ crosses $\Gamma_1$ non tangentially and
 $\Gamma_1 \subset \partial D\times (0,T')$\ .}
\end{assum}
The following result is a particular case of Theorem 4.2 of \cite{FleJ:92}.
\begin{theorem}\label{FJ}
Let $D$ be a bounded open set with a smooth boundary.
Let $N\subset D\times [0,T']\, , T'<T$ satisfy Assumption \ref{(A)}. For the SDE \paref{be}
assume that $\sigma$ is infinitely many times differentiable and bounded,
 the drift $\widetilde b_\ep$ is Lipschitz continuous uniformly with respect
to $\ep$ and enjoys the development
\begin{equation}\label{dev-drift}
\widetilde b_\ep=\widetilde b+\ep \widetilde b_1+o(\ep)\ ,
\end{equation}
uniformly on compact sets of $N$ where $\widetilde b$, $\widetilde b_1$ are $C^\infty$ functions.
Then for $(x,s)\in N$ the following expansion holds
\begin{equation}\label{dev-pr}
\P^\ep_{x,s}(\tau\le T)=\e^{-w(x,s)}\,\e^{-u(x,s)/\ep}\bigl(
1+o(\ep)\bigr)
\end{equation}
uniformly on compact subsets of $N$, where $w: N\to \RR^+$ is the solution of
\begin{equation}\label{sys1}
\begin{cases}
\displaystyle{\partial w\over\partial s}+
\langle \widetilde b-a\nabla u,\nabla w\rangle=-\frac 12\,\tr(a\cdot{\rm Hess}\, u)-
\langle \widetilde b_1,\nabla u\rangle&{\rm in\ }N\cr
w=0&{\rm on\ }\partial D\times [0,T[\cap \overline N\ .
\end{cases}
\end{equation}
\end{theorem}

\begin{remark}\rm The original result in \cite{FleJ:92}
deals with a more general situation in particular providing a full development of $\ep\mapsto\P^\ep_{x,s}(\tau\le T)$.
Beware of some notation mismatch between Theorem \ref{FJ} and the original Theorem 4.2 of
\cite{FleJ:92}; in particular our $\widetilde b_1$ corresponds with $b_2$ there.
\end{remark}
The hypotheses in Theorem \ref{FJ} ensure that for $(x,s)\in N$, there exists a unique
minimizing path for the quantity in the right-hand side of \paref{eq-for-u},
which moreover coincides with the solution $\gamma_{x,s}$ of  the ordinary equation \paref{e159}
for $t\in [s, t^*_{x,s}]$, $t^*_{x,s}$ turning out to be the  first time at which $\gamma_{x,s}$ reaches $\partial D$.
Furthermore, the differential systems \paref{sys1} for $w$ can be solved by characteristics: one has to solve the ordinary
equation \paref{e159} and then
\begin{equation}\label{140}
w(x,s)=\int_s^{t^*_{x,s}} \Bigl(\frac 12\,\tr(a\cdot{\rm Hess}\, u)(\gamma_{x,s}(t),t))+
\langle \widetilde b_1(\gamma_{x,s}(t),t),\nabla u(\gamma_{x,s}(t),t)\rangle\Bigr)\, dt\ .
\end{equation}
\begin{remark}\label{relevant}\rm It is useful to point out two features of Theorem \ref{FJ}.
First, because of Assumption \ref{(A)},
it requires that the characteristic $\gamma_{x,s}$ reaches the boundary at a time
$t^*_{x,s}<T$. 

Second, remark that, by the theory of Large Deviations, the asymptotics as $\ep \to 0$
 of the quantity of interest $\P^\ep_{x,s}(\tau\le T)$ are of the form
$c(\ep)\,\e^{-\ell/\ep}$, where $c$ is a subexponential function of $\ep$, i.e. such that $\lim_{\ep\to0}\ep\log c(\ep)=0$.
Theorem \ref{FJ} states that, under the assumptions considered, the term before the exponential, $c$,
is a constant as a function of $\ep$.

A typical situation when Theorem \ref{FJ} does not apply, for instance,
is when $dX^\ep_t=\sqrt{\ep}\,dB_t$ and $D=]-\infty,L[$ for some $L>0$.
In this case $\gamma_{x,s}(t)=x+\frac{t-s}{1-s}(L-x)$
(so that $\gamma_{x,s}$ reaches the boundary
$\pa D=\{L\}$ at time $T=1$ and the assumptions of Theorem \ref{FJ} are not satisfied) and,
easily by the reflection principle,
$$
\dlines{
\P^\ep_{x,s}(\tau\le 1)=\P\Bigl(\sup_{t\le 1}\  x + \sqrt{\ep} (B_t-B_s) \ge L\Bigr)=\cr
=2\P\Bigl(B_1-B_s\ge \frac{L-x}{\sqrt{\ep}}\Bigr)\sim \frac 2{\sqrt{2\pi(1-s)}L}\,\sqrt{\ep}\,
\e^{-(L-x)^2/2(1-s)\ep}\ ,\cr
}
$$
so that here the term before the exponential is not a constant.

On the other hand the assumptions of Theorem \ref{FJ} are satisfied in
most cases where $X^\ep$ is the time changed bridge of a Diffusion conditioned
to be at some point $y\in D$ at time $\ep$,
in the sense that a ``large'' subset $N$ of $D\times [0,1[$
satisfies Assumption \ref{(A)}.
\end{remark}
\section{Applications and remarks}
In this section we see how to adapt Theorem \ref{FJ} to the case of the asymptotics \paref{dev1} for the exit probability of a
conditioned Diffusion.

A first problem arises from the fact that the drift of the time changed conditioned Diffusion
has a singularity at time $t=1$ (think of the case of the Brownian bridge where
$\widetilde b(x,t)=-\frac x{1-t}$) so that Theorem \ref{FJ} cannot be applied with $T=1$.
This difficulty is easily overcome, as remarked in
\cite{Bal:95}, because it turns out that
\begin{equation}\label{asym}
\widehat\P^{y,t}_{x,s}(\tau<1)\sim \widehat\P^{y,t}_{x,s}(\tau<1-\delta)
\end{equation}
for some $\delta>0$, in the sense that the difference between these two probabilities is
exponentially negligible as $t\to0$. In order to see this, recall
Large Deviation estimates recently obtained for conditioned Diffusions
(see \cite{inahama}, \cite{bailleul} for the case of a compact manifold).
These state that, as $t\to 0$, the time changed conditioned Diffusion starting at $x$ at time $s$
satisfies
a Large Deviation Principle with rate function given by
\begin{equation}\label{Js}
\widehat J_s(\gamma)=\begin{cases}
J_s(\gamma)-J_s(\gamma_0)&\mbox{if }\gamma(1)=y\cr
+\infty&\mbox{otherwise}
\end{cases}
\end{equation}
where
$$
J_s(\gamma)=\frac 12\int_s^1\langle a^{-1}(\gamma(t))\dot\gamma(t),\dot\gamma(t)\rangle\, dt
$$
and $\gamma_0$ denotes a minimizing geodesic (see below) joining $x$ to $y$. Therefore we have
$$
t\log\widehat\P^{y,t}_{x,s}(\tau<1)\sim -\inf_{\gamma_s=x,\tau(\gamma)<1} \widehat J_s(\gamma)\ .
$$
Assume that there exists a unique $\widehat\gamma$ minimizing the right-hand side above, then we can split
\begin{equation}\label{rightmost}
\widehat\P^{y,t}_{x,s}(\tau<1)=\widehat\P^{y,t}_{x,s}(\tau<1,U(\eta,\widehat\gamma))+
\widehat\P^{y,t}_{x,s}(\tau<1,U(\eta,\widehat\gamma)^c)\ ,
\end{equation}
where $U(\eta,\widehat\gamma)$ denotes a neighborhood of radius $\eta$ of the minimizer $\widehat\gamma$.
As the infimum of $\widehat J_s$ on the set of paths $\{\tau<1,U(\eta,\widehat\gamma)^c\}$ is strictly
larger than the
infimum over $\{\tau<1,U(\eta,\widehat\gamma)\}$, the rightmost term in \paref{rightmost} is
exponentially negligible. Let us choose $\eta=\frac 14\, dist(y,\pa D)$ and let $\delta$ be such
that $dist(\widehat\gamma_t, y)\le \eta$ for every $t\ge 1-\delta$. Then for every
$\gamma\in U(\eta,\widehat\gamma)$ and  $t\ge 1-\delta$ we have
$dist(\gamma_t, y)\le 2\eta=\frac 12\, dist(y,\pa D)$. Therefore if $\tau(\gamma)<1$,
then necessarily  $\tau(\gamma)<1-\delta$.

A second question in order to apply the Fleming-James Theorem \ref{FJ} to our problem is the determination of the development of the drift
of the equation for $\eta^t$ in \paref{time-changed}, i.e. of finding vector fields $\widetilde b$ and $\widetilde b_1$ (of course depending on the conditioning point $y$) such that
$$
\bigl(b(\eta^t_v)+\widehat b^{y,t}(\eta^t_v,tv)\bigr)t= \widetilde b(z,v)+t \widetilde b_1(z,v)+o(t),\quad \mbox{as }t\to 0
$$
uniformly on compact sets and then to compute the corresponding quantities $u$ and $w$ of Theorem \ref{FJ}. This in particular requires
to obtain the development of the quantity $\nabla_z \log p(t(1-v),z,y)$, as explained in \S\ref{CD}.

The tool to this goal is provided by Mol{\v{c}}anov results \cite{molchanov-MR0413289} (see also \cite{AZ:81}, Theorem 1.1 p.56) together with those of
Ben Arous \cite{benarous}.
Let us assume that $a=\sigma\sigma^*$ is elliptic. One can then consider on $\R^d$ the Riemannian metric associated to the matrix field $a^{-1}$.
This allows to define the length of a smooth curve $\zeta:[0,1]\to\R^d$ by
$$
l(\zeta)=\int_0^1\sqrt{\langle a^{-1}(\zeta(t))\dot\zeta(t),\dot\zeta(t)\rangle}\, dt
$$
and the corresponding Riemannian distance by
\begin{equation}\label{d1}
d(x,y)=\inf_{\zeta,\zeta(0)=x,\zeta(1)=y}l(\zeta)\ .
\end{equation}
Under an assumption of closeness of the points $x,y$, to be made precise below,
we have (\cite{molchanov-MR0413289}) the
development as $t\to0$
\begin{equation}\label{p56-az}
\log p(t,x,y)\sim-\frac d2\,\log (2\pi t)+\log H(x,y)-\frac 1{2t} d(x,y)^2+A(x,y)
\end{equation}
where $d$ denotes the Riemannian distance \paref{d1} of the metric $a^{-1}$, $H(x,y)=(\det\exp_x'(\xi))^{-1/2}$, $\exp_x$ denoting the exponential map of the
Riemann structure associated to the metric $a^{-1}$ and $\xi$ the tangent vector at $t=0$ of the minimizing geodesic joining $x$ to $y$, and
\begin{equation}\label{only-b}
A(x,y)=\int_0^1\langle a^{-1}(\gamma(t))b(\gamma(t)),\dot\gamma(t)\rangle\, dt
\end{equation}
$\gamma$ denoting again the unique geodesic joining $x$ to $y$. These results are obtained under some regularity assumptions on the coefficients $b$ and $\sigma$ that should be 4 times differentiable.

As mentioned above this development holds under the assumption for the two points $x,y$ to
be close i.e. that they are joined by a unique geodesic along which they are not conjugated. It is a well known fact in Riemannian geometry that for every $x$ there exists a neighborhood $\cl U_x$ of $x$ such that this assumption is satisfied for every $y\in \cl U_x$.

Both $H$ and $d$ are quantities only depending on
the metric $a^{-1}$ and not on the drift $b$ which appears only in the quantity  \paref{only-b}.

Moreover Th\'eor\`eme 3.4 in \cite{benarous} allows to state that the behavior
as $t\to0$ of $\nabla_x \log p(t,x,y)$ is obtained by taking formally the derivatives of the right-hand side in \paref{p56-az}. We have therefore, as $t\to0$,
\begin{equation}\label{dev-der}
t\,\nabla_x \log p(t(1-v),x,y)\sim -\frac 1{2(1-v)}\nabla_x d(x,y)^2 +t(\nabla_x \log H(x,y)+\nabla_x A(x,y))\ .
\end{equation}
We plan, in a forthcoming paper, to use this development in order to be able to obtain explicitly the values
of the constants $c=\e^{-w(x,s)}$,
$\ell=u(x,s)$ appearing in the asymptotics \paref{dev1} for the most common models of stochastic volatility.
In this note, as pointed out at
the beginning, we just wish to investigate the question whether the drift $b$ has an influence
in the development \paref{dev1}. We already know that
the answer is no for a large class of Diffusions in dimension $d=1$ (\cite{BC2002-MR1925452}) and also, in a
multidimensional setting, if $X$ is a Brownian motion with a constant drift:
the bridge of a Brownian motion with a constant drift is exactly equal to a Brownian bridge,
so here to the (constant) drift has no effect.
We start first with an example.
\begin{ex}\label{OU-d}\rm Let
\begin{equation}\label{eds ou}
dX_t=M X_t\, dt+dB_t
\end{equation}
be a $d$-dimensional Ornstein-Uhlenbeck process where $M$ is a $d\times d$-dimensional matrix.
Let us start computing the development of the drift of
the conditioned Diffusion with $X_t=y$ as $t\to0$.

The transition function $p(t,x,\cdot)$ is the density of a $N(\e^{Mt} x,S_t)$-distributed r.v., where
$$
S_t=\int_0^t\e^{Mu}\e^{M^*u}\, du\ .
$$
Therefore
$$
\log p(t(1-v),z,y)=-\frac d2\,\log 2\pi-\frac 12\,\log \det(S_{t(1-v)})-\frac 1{2}\,\langle S^{-1}_{t(1-v)}(y-\e^{Mt(1-v)}z),(y-\e^{Mt(1-v)}z)\rangle
$$
and
$$
\nabla_{z}\,\log p(t(1-v),z,y)=\frac 12\, \e^{M^* t(1-v)}S^{-1}_{t(1-v)}\,(y-\e^{Mt(1-v)}z)+\frac 12\, (y-\e^{Mt(1-v)}z)^*S^{-1}_{t(1-v)}\e^{M t(1-v)}\ .
$$
Writing down the developments as $t\to 0$ of the various terms appearing above we have
$$
\dlines{
S_{t}=tI+(M+M^*)\frac {t^2}2+o(t^2)\ ,\quad
S_{t}^{-1}=\frac 1t\,\bigl(I-(M+M^*)\tfrac {t}2+o(t^2)\bigr)\cr
\e^{tM^*}=I+tM^*+o(t),\quad \e^{tM}=I+tM+o(t)\cr
}
$$
so that
\begin{align*}
&\e^{tM^*}S_{t}^{-1}=\frac 1t\,I-\frac 12\,(M+M^*)+M^*+o(1)=\frac 1t\,I+\frac 12\,(M^*-M)+o(1)\cr
&S_{t}^{-1}\e^{tM}=\frac 1t\,I-\frac 12\,(M+M^*)+M+o(1)=\frac 1t\,I+\frac 12\,(M-M^*)+o(1)\ .\cr
\end{align*}
Also
$y-\e^{M t}z=y-z+z-\e^{M t}z=y-z-Mtz+o(t)$, hence
$$
\dlines{
\frac 12\, \e^{M^* t(1-u)}S^{-1}_{t(1-u)}\,(y-\e^{Mt(1-u)}z)
=\frac 12\, \Bigl(\frac 1{t(1-u)}\,I+\frac 12\,(M-M^*)+o(1)\Bigr)(y-z-Mt(1-u)z+o(t))=\cr
=\frac 12\,\Bigl(\frac {y-z}{t(1-u)}+\frac 12\,(M-M^*)(y-z)-Mz+o(1)\Bigr)\ .
}
$$
Similarly
$$
\dlines{
\frac 12\, (y-\e^{Mt(1-u)}z)^*S^{-1}_{t(1-u)}\e^{M t(1-u)}
=\frac 12\,\Bigl(\frac {y-z}{t(1-u)}+\frac 12\,(M^*-M)(y-z)-M^*z+o(1)\Bigr)
}
$$
and putting things together, after some straightforward computations we find
$$
\dlines{
t\nabla_{z}\,\log p(t(1-v),z,y)
=\frac {y-z}{1-v}-\frac 12\,t(M+M^*)z+o(t)\ .
}
$$
Remark that the same result would have been obtained very quickly also computing the terms appearing in \paref{dev-der}, as here
$H(x,y)\equiv 1$ and $d(x,y)=|x-y|$. Therefore the asymptotics for the drift of the bridge of $X$ given $X_t=y$ is
$$
t\bigl(b(z)+\widehat b^{y,t}(z,tv)\bigr)=\frac {y-z}{1-v}-\frac 12\,t(M+M^*)z+tM z+o(t)=\frac {y-z}{1-v}+\frac 12\,t(M-M^*)z+o(t)\ .
$$
Hence we are in the situation of \paref{dev-drift} with
$$
\widetilde b(z,v)=\frac {y-z}{1-v},\quad \widetilde b_1(z,v)=\frac 12\,(M-M^*)z\ .
$$
Remark that $\widetilde b_1\equiv0$ if and only if the matrix $M$ is symmetric. Therefore, in general, the quantity $w$ which determines the value
of the constant $c$ in the expansion \paref{dev1} will depend on the drift $z\mapsto Mz$, if the matrix $M$ is not symmetric.

To be more precise let us consider the case where $D$ is the half-space
$\{z,\langle \vec{v},z\rangle < k\}$ for some $\vec{v}\in\R^d$, $|\vec{v}|=1$ and $k>0$.
Let $x,y\in  D$ and let us evaluate the expansion of
\paref{dev1} for the process $X$ conditioned by $X_t=y$, where $\tau$ denotes the exit time from $D$.

Note first that the function $u$ defined in \paref{eq-for-u} coincides with the one
for the bridge of the Brownian motion, i.e.
\begin{equation}\label{bridge u}
u(x,s)=\frac{2}{1-s} \left(k-\langle x,\vec{v} \rangle \right) \left(k-\langle y,\vec{v} \rangle \right)
\end{equation}
(moreover it is easy to check that such a function $u$ satisfies \paref{HJ1}). Of course
we have  $\Delta u \equiv 0$ and
$$
\nabla u(x,s) = -\frac{2}{1-s}\left(k-\langle y,\vec{v} \rangle\right)\vec{v}\ .
$$
Therefore the sharp asymptotics for the exit time probability $q_t(x,s)$ of \paref{dev1},
as the conditioning time $t$ tends to $0$, for the Diffusion starting at $X_s=x$ and
conditioned by $X_t=y$ is
\begin{equation}
q_t(x,s) \sim \e^{-w(x,s)}\e^{ -\frac{u(x,s)}{t}}\ ,
\end{equation}
where, recalling \paref{140},
\begin{equation}\label{w1}
w(x,s) = -\left(k-\langle y,\vec{v} \rangle\right)\int_s^\tau \frac{1}{1-t}
\langle (M-M^*)\gamma_{x,s}(t), \vec{v} \rangle\,dt\ ,
\end{equation}
$\gamma_{x,s}$ being the solution of
\begin{equation}\label{eq gamma}
\dot \gamma_{x,s}(t) =\widetilde b(\gamma_{x,s}(t),t) -\nabla u(\gamma_{x,s}(t),t)\ ,\qquad \gamma_{x,s}(s)=x
\end{equation}
and $\tau=t^*_{x,s}$ the time at which $\gamma_{x,s}$ reaches the boundary $\pa D$.
Straightforward computations lead to the solution
\begin{equation}
\gamma_{x,s}(t) = x+\frac{t-s}{\tau-s}(\eta-x)\,\qquad s\le t \le \tau\ ,
\end{equation}
where
\begin{equation}\label{ttt}
\tau = s+(1-s)\frac{k-\langle x,\vec{v}\rangle}{2k - \langle  x+y  , \vec{v}  \rangle}\quad\mbox{and}\quad \eta =
x+ \frac{k-\langle x,\vec v\rangle}{2k - \langle x+y, \vec v\rangle}
(y-x+2(k-\langle y, \vec{v} \rangle)\vec{v})\ .
\end{equation}
Remark that  $\tau<1$, $\eta\in \pa D$ and does not depend on $s$, furthermore $\gamma_{x,s}$ is the
line segment connecting $x$ with $\eta$.
Going back to \paref{w1} we have
$$
w(x,s)=-\left(k-\langle y,\vec{v} \rangle  \right)\int_s^\tau \frac{1}{1-t}\langle \vec{v}, (M-M^*)\left(
x + \tfrac{t-s}{\tau-s}(\eta -x)
\right)
   \rangle\,dt
$$
which gives easily
\begin{equation}\label{w-explicit}
 w(x,s)
=\bigl(k-\langle y,\vec{v} \rangle\bigr)\Bigl[\frac{k-\langle x,\vec{v}\rangle}{2k - \langle  x+y  , \vec{v}  \rangle}\langle \vec{v}, (M-M^*)(y -x)   \rangle+\log \Bigl(\frac{k-\langle y,\vec{v}\rangle}{2k - \langle  x+y  , \vec{v}  \rangle} \Bigr)\,\langle \vec{v}, (M-M^*)y
   \rangle  \Bigr]
\end{equation}
($w(x,s)$ does not depend on $s$).
Therefore for the quantity $q_t(x,s)$ in \paref{dev1} we have
\begin{align*}
&c_{x,s}=\Bigl(\frac{k-\langle y,\vec{v}\rangle}{2k - \langle  x+y  , \vec{v}  \rangle}\Bigr)^{-(k-\langle y,\vec{v} \rangle)\langle \vec{v}, (M-M^*)y
   \rangle}\,\e^{-\frac{(k-\langle y,\vec{v} \rangle)(k-\langle y,\vec{v} \rangle)}{2k-\langle x+y,v\rangle}\,\langle \vec{v}, (M-M^*)(y -x)
   \rangle}\cr
&\ell_{x,s}=\frac 2{1-s} \left(k-\langle x,\vec{v} \rangle\right) \left(k-\langle y,\vec{v} \rangle\right)\ .
\end{align*}
We stress that the expansion \paref{dev1}  depends indeed on the matrix $M$, if this is not symmetric. Whereas if $M$
is any symmetric matrix, its value has
no influence on the expansion \paref{dev1}, which is then exactly the same as if $X$ was the Brownian motion, where $w\equiv0$.

We did not bother to check the assumptions of Theorem \ref{FJ}. It is not however difficult, given the computations above, for a given $x\in D$,
to construct a RSR containing $(x,0)$. Indeed remark that the expression of $\tau$ in \paref{ttt} says that every characteristic
$\gamma_{x,s}$ reaches $\pa D$ at a time that is strictly smaller than $1$. One can therefore construct a RSR of the form
$N=\{(z,s);z\in D, t^*_{z,s}<T'\}$, for some $T'$ such that $t^*_{x,0}<T'<1-\delta$ where $\delta$ is given in \paref{asym}. The only remaining assumption to be checked is that $D$ is assumed there to be
bounded, which is not our case. This point is explained in the next remark.
\end{ex}
\begin{remark}\label{rem bounded}\rm
It is easy to show that the asymptotics for the probability of exit from an open set $D$ is, by a standard localization argument, the same
as for the exit from a suitable bounded subset $\widetilde D \subset D$.

To be more precise, a repetition of the argument leading to \paref{asym} yields
$$
q_t(x,0)\sim \widehat \P^{y,t}_{x,0}(\tau < 1, U(\eta, \widehat \gamma))\ ,
$$
where $U(\eta, \widehat \gamma)$ is a neighborhood of radius $\eta$ of the minimizer $\widehat \gamma$
for $\inf_{\gamma_0=x, \tau(\gamma)<1} \widehat J_0(\gamma)$, $\widehat J_s$ being
given in \paref{Js}. One can then set $\widetilde D$ to be the intersection of $D$ with a bounded neighborhood
of the support of $\widehat\gamma$, chosen in such a way as to preserve the smoothness of the boundary.
\end{remark}
The previous example shows, among other things, that the sharp Large Deviation estimate of exit of the
bridge of a multidimensional
Ornstein-Uhlenbeck process depends on the drift of the original process, unless the matrix $M$ is
symmetric. This is a phenomenon that is
better investigated in the following statement.
\begin{prop}
Let $X$ be the $d$-dimensional diffusion process that is the solution of the SDE
\begin{equation}\label{sde X}
dX_t=b(X_t)\, dt+\sigma(X_t)\, dB_t
\end{equation}
and assume that $a=\sigma\sigma^*$ is elliptic and that $b$ and $\sigma$ are $4$ times differentiable.
Let us denote $\eta^t$ the corresponding process conditioned given $X_t=y$, $t>0$ and time changed (see \paref{time-changed}).
Let $x$ be close enough to $y$, in the sense that $x$ and $y$ are joined by unique geodesic $\gamma_0$ of the metric $a^{-1}$
along which they are not conjugated.
Then if there exists a potential $U:\R^d \to \R$ such that $\nabla\, U=a^{-1}b$,
the development for the drift of $\eta^t$ up to the first order as $t\to0$ does not depend on $b$.
\end{prop}
\noindent {\it Proof.} We must prove for the drift $(b(x)+\widehat b^{y,t}(x,tv))t$ of the time changed conditioned process $\eta^t$ a
development of the form $(b(x)+\widehat b^{y,t}(x,tv))t=\widetilde b+t\widetilde b_1+o(t)$ where neither $\widetilde b$ nor
$\widetilde b_1$ depend on $b$. Recall the development \paref{dev-der}: thanks to \paref{only-b}, if $\nabla\, U=a^{-1}b$ we have of course
$$
A(x,y)=U(\gamma_1)-U(\gamma_0)=U(y)-U(x)
$$
and $\nabla_x A(x,y)=-\nabla U(x)=-a^{-1}(x)b(x)$. Hence, by \paref{dev-der}, as $t\to0$,
$$
\widehat b^{y,t}(x,tv)=a(x)\nabla_x\log p(t(1-v),x,y)\sim a(x)\Bigl(\nabla_x\log H(x,y)-\frac 1{2t(1-v)}\nabla_x d(x,y)^2- a^{-1}b(x)\Bigr)
$$
so that the drift of the time changed conditioned process is
$$
\dlines{
\bigl(b(x)+\widehat b^{y,t}(x,tv)\bigr)t=tb(x)+ta(x)\Bigl({\nabla_x} \log H(x,y)-
\frac 1{2t(1-v)}\nabla_x  d(x,y)^2- a^{-1}b(x)\Bigr)=\cr
=t\nabla_x \log H(x,y)-\frac 1{2(1-v)}\nabla_x  d(x,y)^2\ ,\cr
}
$$
thus $\widetilde b(x,v)=-\frac 1{2(1-v)}\nabla_x  d(x,y)^2,$ whereas $\widetilde b_1(x,v)= \nabla_x \log H(x,y)$ and neither of these terms
depends on $b$.
\qed
\begin{remark}\rm
If the  hypotheses in Theorem \ref{FJ} are satisfied, then the drift of the unconditioned Diffusion does not influence the sharp
asymptotics for the  probability of exit from the domain $D$, i.e. neither of the functions $u$ and $w$ in \paref{dev-pr} depend on $b$ in \paref{sde X}.
Recall that $D$ could be unbounded, actually the argument in Example \ref{OU-d} holds in great
generality.
\end{remark}
Coming back to Example \ref{OU-d}, of course if $M$ is symmetric then the drift $z\mapsto Mz$ turns out to be the gradient field of the potential $U(z)=\frac 12\,\langle Mz,z\rangle$.
We have therefore proved that, whenever the Fleming-James Theorem~\ref{FJ} can be applied, the asymptotics \paref{dev1} do not depend on the drift $b$ of the original Diffusion as far
as $a^{-1}b$ is a gradient field and also (Example \ref{OU-d}) that the drift can be influential if this assumption is not satisfied.
\def\cprime{$'$} \def\cprime{$'$}

\end{document}